\documentclass[a4paper,12pt]{article}
\usepackage{graphicx, amssymb, amsmath, amsthm, mathrsfs, fullpage, hyperref, tikz, pgf, epstopdf, enumerate, array, cite, booktabs, lscape, pdflscape,multirow,makecell,enumitem}
\hypersetup{
    plainpages=true,
    pdfstartview=FitH,
    bookmarksopen=true,
    colorlinks=true,
    linkcolor=blue,
    citecolor=blue
}

\newtheorem{thm}{Theorem}[section]

\newtheorem{lmm}[thm]{Lemma}

\newtheorem{pro}[thm]{Proposition}

\newtheorem{rem}[thm]{Remark}

\newcommand*{\myand}{\,\and\,}

\begin{document}

\title{\textbf{Asymptotic Enumeration of\\ Subclasses of Level-$2$ Phylogenetic Networks 
}}
\author{Hexuan Liu\thanks{Department of Pure and Applied Mathematics, Graduate School of Fundamental Science and Engineering, Waseda University, 169-8555 Tokyo, Japan. Email: \href{lhx@ruri.waseda.jp}{lhx@ruri.waseda.jp}} 
\myand
Bing-Ze Lu\thanks{Department of Applied Mathematics, National Sun Yat-sen University, Kaohsiung 804, Taiwan. Email: \href{wind951016@gmail.com}{wind951016@gmail.com}}
\myand
Taoyang Wu\thanks{School of Computing Sciences, University of East Anglia (UEA), Norwich, NR4 7TJ, UK. Email: \href{taoyang.wu@uea.ac.uk}{taoyang.wu@uea.ac.uk}}
\myand
Guan-Ru Yu\thanks{Department of Applied Mathematics, National Sun Yat-sen University, Kaohsiung 804, Taiwan. Email: \href{gryu@math.nsysu.edu.tw}{gryu@math.nsysu.edu.tw}}}
\date{\today}
 \maketitle

\begin{abstract}
This paper studies the enumeration of seven subclasses of level-$2$ phylogenetic networks under various planarity and structural constraints, including terminal planar, tree-child, and galled networks. We derive their exponential generating functions, recurrence relations, and asymptotic formulas. Specifically, we show that the number of networks of size $n$ in each class follows:
\[
N_n \sim c \cdot n^{n-1} \cdot \gamma^n,
\]
where $c$ is a class-specific constant and $\gamma$ is the corresponding growth rate. 
Our results reveal that being terminal planar can significantly reduce the growth rate of general level-2 networks, but has only a minor effect on the growth rates of tree-child and galled level-2 networks. Notably, the growth rate of 3.83 for level-$2$ terminal planar galled tree-child networks is remarkably close to the rate of 2.94 for level-$1$ networks.

\end{abstract}
\section{Introduction}
Phylogenetic trees have been instrumental in evolutionary biology, illustrating how species diverge from common ancestors. However, their limitations in representing events like hybridization and horizontal gene transfer led to the development of phylogenetic networks. These networks extend trees by incorporating complex evolutionary events, offering a more flexible framework for modeling evolutionary histories.

A \textbf{phylogenetic network} with $n$ leaves is an acyclic directed graph with a single source (the root) and a set of sinks (leaves) bijectively labeled by $\{1, 2, \dots, n\}$. Throughout this paper, all phylogenetic networks are binary, satisfying the following conditions:
\setlist[itemize]{itemsep=0.5pt, topsep=5pt} 
\begin{itemize}
    \item[(i)] If $n > 1$, the root has exactly one child, and a leave has exactly one parent.
    \item[(ii)] Each interior node is either a tree node (one parent and two children) or a reticulation node (two parents and one child).
    \item[(iii)] For every 2-connected component $B$ with three or more vertices, there are at least two cut arcs directed away from $B$.
\end{itemize}

A phylogenetic network is a \textbf{tree-child network} if every interior node has at least one child that is either a tree node or a leaf (see \cite{CaZh,FuYuZh1}). Similarly, it is a \textbf{galled network} if every reticulation node lies in a unique tree cycle (see \cite{FuYuZh2,GuRaZh}). A \textbf{galled tree-child network} (GTC) is a phylogenetic network that is both a galled network and a tree-child network (see\cite{ChFuYu}).

\begin{figure}[h!]
    \centering
    \includegraphics[width=0.6\textwidth]{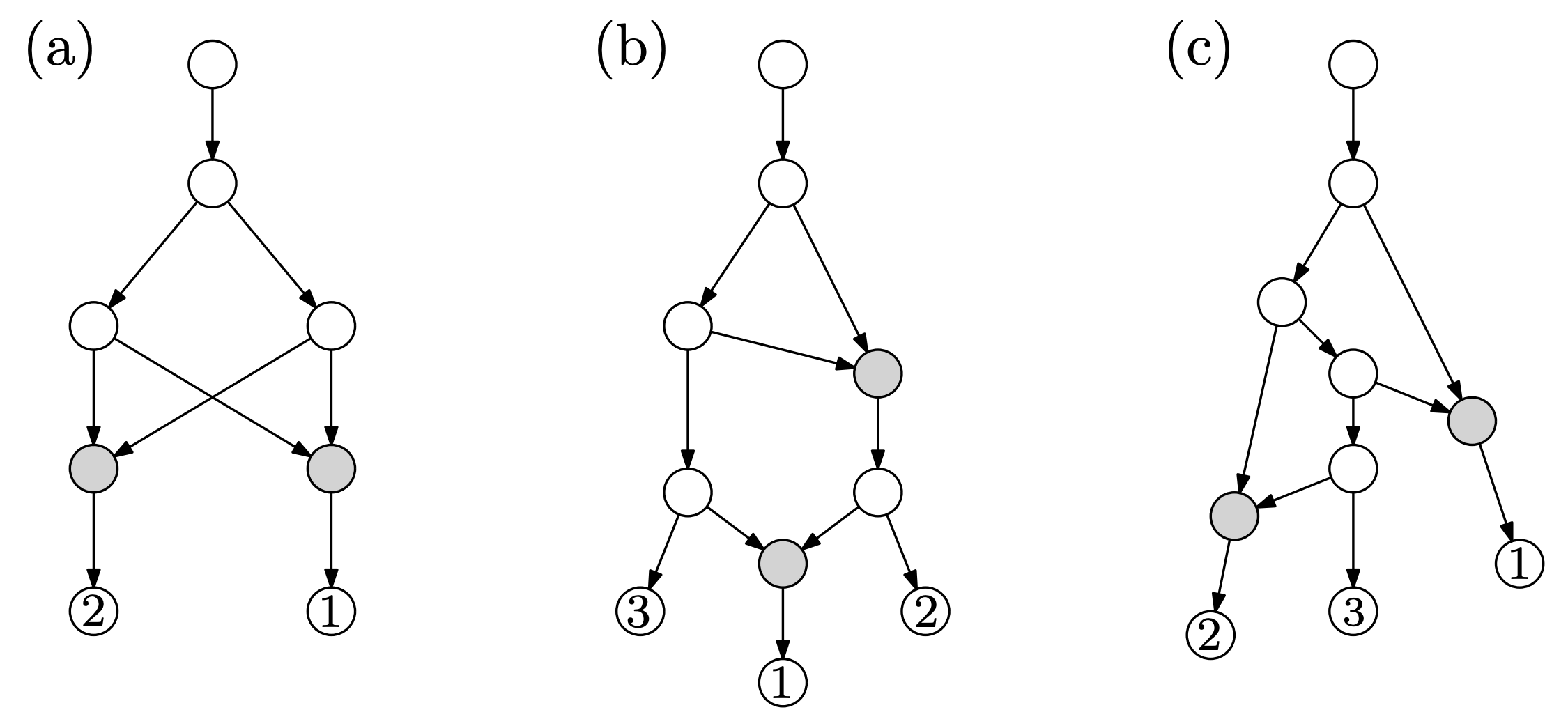}
        \vspace{-10pt} 
    \caption{Examples of level-2 networks where reticulation nodes are represented by filled circles. (a) is a galled network which is not tree-child; (b) is a tree-child network which is not galled; (c) is a galled tree-child network. 
    Both (b) and (c) are terminal and outer planar, whereas (a) is neither.}
    \label{fig:gtc}
\end{figure}

A particularly important class is \textbf{level-$k$ networks}, where each 2-connected component contains at most $k$ reticulation nodes (see \cite{GaBePa}). See Figure~\ref{fig:gtc} for examples of level-2 networks. This hierarchical classification controls the complexity of networks, making them computationally tractable. For instance, level-$1$ networks correspond to galled trees (see \cite{HuIeMoSc,HuMoSeWu}), while higher levels (such as level-$2$) model more complex evolutionary histories (see\cite{BoGaMa,IeKeKeStHaBo,IeMo}).

The \textbf{planarity} of phylogenetic networks is crucial for their visualization and interpretation. A phylogenetic network is \emph{planar} if it can be embedded on an Euclidean xy-plane without edge crossings, except at their endpoints. Based on planarity, phylogenetic networks are classified into four types:
\setlist[itemize]{itemsep=0.5pt, topsep=5pt} 
\begin{itemize}
    \item \textbf{Planar networks} ($\mathcal{N}_p$): Admit a planar embedding with no edge crossings.
    \item \textbf{Upward planar networks} ($\mathcal{N}_u$): 
    Networks with a planar drawing in which each directed  edge 
    is a curve of monotonically non-decreasing in the y-coordinate.
    \item \textbf{Terminal planar networks} ($\mathcal{N}_t$): Planar networks with a planar drawing in which all leaves and the root are on the outer face.
    \item \textbf{Outer planar networks} ($\mathcal{N}_o$): Planar networks with a planar drawing in which all vertices are on the outer face boundary.
\end{itemize}
These classes satisfy the inclusion relation: $\mathcal{N}_p \supset \mathcal{N}_u \supset \mathcal{N}_t \supset \mathcal{N}_o$.

In fact, level-$1$ or $2$ constraints are closely tied to planarity, leading to the following properties (see \cite{MoWu}).

\begin{pro}\label{prop:combined}
The following statements hold:
\setlist[itemize]{itemsep=0.5pt, topsep=5pt} 
\begin{itemize}
\item[(i)] All level-1 networks are outer planar.
\item[(ii)] All level-2 networks are upward planar. Moreover, a level-2 network N is terminal planar if and only if N is outer planar.
\end{itemize}
\end{pro}

Proposition~\ref{prop:combined} establishes that level-2 tree-child networks are not only planar but also upward planar. Under the level-2 constraint, network planarity reduces to two cases: outer planar networks and non-outer planar networks (which remain upward planar).

It is well known that there are $(2n-3)!!$ binary phylogenetic trees with $n$ leaves. However, counting phylogenetic networks is considerably more challenging. Consequently, much of the recent progress has focused on subclasses with certain topological constraints, such as level-1 and level-2 networks~\cite{BoGaMa}, tree-child networks~\cite{mcdiarmid2015counting,FuYuZh1}, and galled networks~\cite{GuRaZh,FuYuZh2}.

This paper focuses on enumerating subclasses of level-2 networks, including \textbf{level-2 tree-child networks} and \textbf{level-2 galled networks} under various planarity constraints. We derive functional equations for their exponential generating functions and asymptotic growth rates, revealing their structural properties and combinatorial behavior. Additionally, we study \textbf{galled tree-child networks}, the intersection of the two classes above. Our results deepen the mathematical understanding of level-$2$ networks. In Section~\ref{sec:main-results}, we present the main results, followed by the methodology in Section~\ref{sec:gf}, where generating functions and theoretical tools are employed for enumeration.

\section{Main Results}\label{sec:main-results}

\begin{thm}\label{thm:main}
Let $t_n$ denote the number of level-2 tree-child networks with $n$ leaves. Then, the asymptotic growth of $t_n$ is given by
\begin{align}
\label{eq:main:growth}
t_n \sim c \cdot n^{n-1} \cdot \gamma^n,
\end{align}
where $c \approx 0.0667418464$ and $\gamma \approx 4.6710490708$. Moreover,  each of the seven subclasses of level-2 networks in the following table has the same asymptotic growth pattern as in Eq.\eqref{eq:main:growth}  with distinct constants $c$ and growth rates $\gamma$ as summarized in the following table. 

\begin{center}
\begin{minipage}{\textwidth} 
\centering
\begin{tabular}{c|c|c|c|c}
\toprule
 & \textbf{general} & \textbf{tree-child} & \textbf{galled} & \textbf{GTC} \\ 
\midrule
\makecell[l]{arbitrary\\
planar\\ upward planar} 
 & \makecell[l]{$c=0.02931010$\\ $\gamma=15.4332995$}
 & \makecell[l]{$c=0.06674185$\\ $\gamma=4.67104907$} 
 & \makecell[l]{$c=0.05885954$\\ $\gamma=6.42241234$} 
 & \makecell[l]{$c=0.07888067$\\ $\gamma=3.98275804$} \\
\midrule
\makecell[l]{terminal planar\\ outer planar} 
 & \makecell[l]{$c=0.03486095$\\ $\gamma=12.9230111$} 
 & \makecell[l]{$c=0.07450612$\\ $\gamma=4.33252428$} 
 & \makecell[l]{$c=0.06965278$\\ $\gamma=5.39994365$} 
 & \makecell[l]{$c=0.08586449$\\ $\gamma=3.83201916$} \\
\bottomrule
\end{tabular}
\end{minipage}
\end{center}
\end{thm}

\begin{rem}
The first cell follows from Bouvel et al.~\cite{BoGaMa}: they showed that the exponential growth rate for level-1 networks is $2.94$ and that for  level-2 networks is $\gamma = 15.43$. In contrast, we found that the rate for level-2 terminal planar networks is $\gamma = 12.92$, indicating a noticeable gap. However, under the tree-child or galled constraints, this gap is significantly smaller. Moreover, the growth rate of level-2 galled tree-child networks ($\gamma = 3.98/3.83$) is rather close to that of level-1 networks. This suggests that the tree-child and galled constraints impose limitations similar to those of level-1 networks. 
\end{rem}

\section{Asymptotics via Generating Functions}\label{sec:gf}
In this section, we describe a general methodology for deriving asymptotic results on the number of networks from their generating functions. We illustrate this approach using level-2 tree-child networks as a primary example to outline a proof for Theorem~\ref{thm:main}. 

\begin{figure}[h!]
    \centering
    \includegraphics[width=0.8\textwidth]{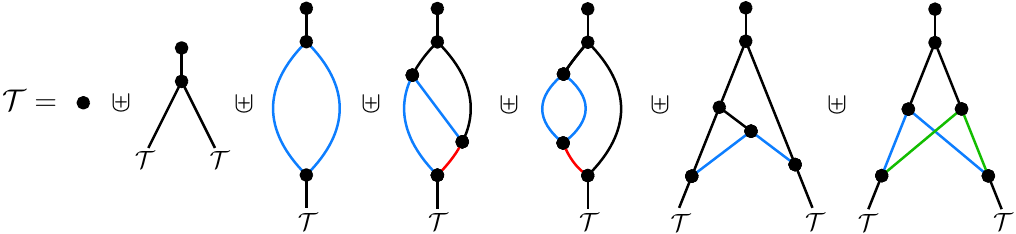}
        \vspace{-10pt} 
    \caption{Decomposition of the combinatorial structure for level-2 tree-child networks. Red edges indicate that they cannot be empty (must contain some insertion). In contrast, pairs of blue or green edges represent that, within each pair of edges of the same color, both cannot be empty simultaneously.}
    \label{fig:breakdown}
\end{figure}

\begin{lmm}\label{lmm:fegf}
Let $t_n$ denote the number of level-$2$ tree-child networks with $n$ leaves, and let $T(x)$ be the corresponding exponential generating function defined as
\begin{align}\label{egf1}
T(x) = \sum_{n \geq 1} t_n \frac{x^n}{n!}.
\end{align}
Then, $T(x)$ satisfies the functional equation:
\begin{align*}
T ={} & x + \frac{T^2}{2} + \frac{1}{2}\left(\left(\frac{1}{1-T}\right)^2 - 1\right)T + \frac{3}{2}\left(\frac{1}{1-T}\right)^2 \left(\frac{T}{1-T}\right) \left(\left(\frac{1}{1-T}\right)^2 - 1\right)T \\
& + \left(\frac{1}{1-T}\right)^4 \left(\left(\frac{1}{1-T}\right)^2 - 1\right)T^2 + \frac{1}{4}\left(\frac{1}{1-T}\right)^2 \left(\left(\frac{1}{1-T}\right)^2 - 1\right)^2 T^2,
\end{align*}
which is derived from the decomposition illustrated in Figure~\ref{fig:breakdown}.
\end{lmm}

Due to space constraints, we omit the detailed derivation of the generating function (only the decomposition is shown), as it is closely tied to the combinatorial structure of the problem and involves a lengthy process. Instead, we present the functional equation obtained from our analysis. By comparing the coefficients on both sides of the equation and iteratively substituting the results from the left-hand side into $T$ on the right-hand side, we derive the following series expansion:
\[
T(x) = x + \frac{3}{2}x^2 + 11x^3 + \frac{745}{8}x^4 + \frac{3333}{4}x^5 + \frac{127395}{16}x^6 + \frac{637635}{8}x^7 + \frac{105501013}{128}x^8 + \cdots.
\]
This expansion corresponds to the EGF (\ref{egf1}), yielding the sequence
\[
\{t_n\}_{n \geq 1} = \{1, 3, 66, 2235, 99900, 5732775, 401710050, 33232819095, \cdots\}.
\]

Next, we introduce the key tools used to analyze the $n$-th order coefficients of the generating function. We begin with the following version of singular inversion theorem.

\begin{thm}\label{thm:singular}{(Singular inversion theorem)\cite{FlSe}}
Let $C(z)$ be a generating function with $C(0) = 0$, satisfying $C(z) = z \phi(C(z))$, where $\phi(z) = \sum_{n \geq 0} \phi_{n} z^n$ is a power series such that: (i) $\phi_{0} \neq 0$, (ii) all $\phi_{n}$ are non-negative reals, (iii) $\phi(z)  \not\equiv \phi_{0} + \phi_{1}z$ (i.e., $\phi$ is nonlinear).

Let $R$ be the radius of convergence of $\phi$ at $0$. Assume (iv) $\phi$ is analytic at $0$ (so $R > 0$), (v) the equation $\phi(z) - z \phi^{\prime}(z) = 0$ has a unique solution $\tau \in (0, R)$, and (vi) $\phi$ is aperiodic.

Then:
\setlist[itemize]{noitemsep, topsep=3pt} 
\begin{itemize}
    \item The radius of convergence of $C(z)$ is $\rho = \frac{\tau}{\phi(\tau)}$.
    \item Near $\rho$, $C(z) \sim \tau - \sqrt{ \frac{2\phi(\tau)}{\phi^{\prime\prime}(\tau)} } \sqrt{1 - \frac{z}{\rho}}$.
    \item $
[z^n] C(z) = \frac{1}{n} [z^{n-1}] \phi(z)^n.
$
    Furthermore, as $n \to \infty$, $[z^n] C(z) \sim \sqrt{ \frac{\phi(\tau)}{2\phi^{\prime\prime}(\tau)} } \frac{ \rho^{-n} }{ \sqrt{\pi n^3} }$.
\end{itemize}
\end{thm}

We now apply Theorem~\ref{thm:singular} to the functional equation derived from Lemma~\ref{lmm:fegf}. First, we rewrite the equation into a form compatible with Theorem~\ref{thm:singular}, namely $T(z) = z \phi(T(z))$. The function $\phi(z)$ is given by:
\begin{small}
\begin{align*}
\phi(z) &= \frac{-4(z-1)^6}{2z^7 - 18z^6 + 67z^5 - 126z^4 + 124z^3 - 70z^2 + 30z - 4} 
&= 1 + \frac{3}{2}z + \frac{35}{4}z^2 + \frac{403}{8}z^3 + \cdots.
\end{align*}
\end{small}
We can verify that $\phi_{0} \neq 0$ and that $\phi$ is nonlinear. In fact, it is straightforward to check that all coefficients $\phi_{n}$ are positive rational numbers. By solving the characteristic equation 
\[
\phi(z) - z \phi^{\prime}(z) = 0,
\]
we obtain a unique solution $\tau \approx 0.1226285445$ (computed using Maple and the Newton method). From this, we further derive
\[
\rho = \frac{\tau}{\phi(\tau)} \approx 0.0787573489.
\]
The asymptotic expression for $t_n$ is then given by
\[
t_n = n! \cdot [x^n] T(x) \sim \sqrt{ \frac{\phi(\tau)}{\phi^{\prime\prime}(\tau)} } \, n^{n-1} (\rho e)^{-n},
\]
where, after calculation, we obtain the constants in 
Theorem~\ref{thm:main}:
\[
\sqrt{ \frac{\phi(\tau)}{\phi^{\prime\prime}(\tau)} } \approx 0.0667418464 \quad \text{and} \quad \frac{1}{\rho e} \approx 4.6710490707.
\]

\section*{Acknowledgement}
We thank the reviewers of the Proceedings of the 13th European Conference on Combinatorics, Graph Theory and Applications (EUROCOMB'25) for accepting this paper as a conference paper in EUROCOMB'25. We also appreciate their valuable suggestions and guidance on the content.

\end{document}